\documentclass[11pt, a4paper]{article}
\usepackage{color,pdfsync}
\usepackage{times,amsmath,amssymb,amsfonts}
\usepackage{url}
\usepackage{a4wide}
\usepackage{epsfig,epic,eepic,graphicx}
\newtheorem{theorem}{Theorem}

\newtheorem{proposition}[theorem]{Proposition}

\makeatletter

\@addtoreset{equation}{section}
\makeatother

\newcommand{\R}{{\mathbb R}}

\newcommand{\pa}{{\partial}}
\newcommand{\na}{{\nabla}}

\newcommand{\eps}{{\varepsilon}}
\newcommand{\tbb}{\alpha_P}
\newcommand{\tbt}{\alpha_S}

\def\div{\hbox{div  }}

\def\rmx{{\rm x}}
\def\rmX{{\rm X}}

\title{Computation of the drag force on a rough sphere close to a wall}
\author{David G\'erard-Varet\footnote{Institut de Math\'ematiques de Jussieu, 175 rue du Chevaleret, 75013 Paris. \newline \indent \indent Mail: \texttt{gerard-varet@math.jussieu.fr}} \hspace{0.5cm} Matthieu Hillairet\footnote{Ceremade, Place du Mar\'echal de Lattre de Tassigny, 75775 Paris Cedex 16.\newline \indent \indent Mail: \texttt{hillairet@ceremade.dauphine.fr}}}

\begin{document}
\maketitle
\begin{abstract}
We consider the  effect of surface roughness on solid-solid contact in a Stokes flow. Various models for the  roughness are considered, and a unified methodology is given to derive the corresponding asymptotics  of the drag force. In this way, we recover and  clarify the various expressions  that can be found in the litterature. 
\end{abstract}

\section{Introduction}
The dynamics of solid particles in a viscous fluid is crucial to many phenomena, such as blood flow, sedimentation or filtration. The  drag force exerted by the fluid on the solids plays of course a central role  in this dynamics.  It has been the matter of many studies. The first ones focused on  the dynamics of a rigid sphere near a plane wall,  that moves in a Stokes flow  under no-slip conditions: we refer to  the pioneering works \cite{Brenner&Cox63,ONeill64,Cooley&Oneill69,ONeill&Stewartson67}. The main conclusion of these works is that the drag force is inversely proportional to the distance $h=h(t)$ between the sphere and the plane at time $t$. The reduced ordinary differential equation that governs the movement of the sphere is then of the type:
{ $ \ddot{h} + {\dot{h}}/{h} = f$}, which prevents collision between the sphere and the wall in finite time. We quote that this striking conclusion holds for any value of the fluid  viscosity and of the sphere density. Moreover, it is still valid for  arbitrary solids with  smooth surfaces, and it is still valid within an unsteady Navier-Stokes flow (see \cite{Hillairet07}). 

\medskip
This theoretical no-collision result, that goes against Archimedes' principle,  is clearly unrealistic at the scale of macroscopic solids.  Even at microscopic scales, "dry collisions" have been clearly recognized. Therefore, many articles have tried to identify the flaw of the previous modelling, in order to circumvent the paradox.  Among possible flaws that have been suggested one can mention: 
\begin{itemize}
\item The rigidity assumption. Elasticity, even weak could allow for solid contact: see \cite{Davis3}. 
\item The no-slip condition, that  is no longer valid  when the distance between the solids is  of the order of the mean free path of the fluid particles:  see \cite{Hocking}.
\item the incompressibility  assumption{ : see} \cite{Barnocky}. 
\end{itemize} 
We shall focus here on another very  popular explanation for the no-collision paradox: roughness. The basic idea is that nothing is as smooth as a plane or a sphere: irregularity of the surface can thus affect the fluid-solid interaction. This {\em credo} has led to many experimental and theoretical studies, focusing on roughness-induced effects on drag forces (\cite{Smart,Davis1,Feuillebois}). Such studies  will be discussed in the core of the paper. 

\medskip
We quote that the interest in  roughness issues has been renewed these last years, notably in connection to microfluidics. Indeed, it has been recognized that the classical no-slip boundary condition, which is relevant at the macroscopic scale, may fail at the micro- or nanoscale. This happens for instance for some corrugated hydrophobic surfaces, which trap gas bubbles in their humps and generate in this way some substantial slip. More generally, to determine the appropriate boundary condition at a rough surface is a matter of current debate. In this context, if one has  theoretical formulas that express how  the drag   force depends on the "rough" boundary conditions, one  may check experi\-men\-tally through the force measurement what  the right boundary condition is. This interesting point of view is  for instance developed in \cite{Vinogradova2, Vinogradova1}. 

\medskip
The aim of this paper is to investigate mathematically and in a unified way the relation between the roughness and the drag force.  Namely, we study  the evolution with time $t$  of a rough solid $S(t)$, falling  towards a rough wall $P$ in a Stokes flow. We assume for simplicity that the solid moves by translating along the vertical axis $r=0$, where  $(r,\theta,z)$ are  cylindrical coordinates.  We shall comment on this simplification later on. 
Various models for the roughness are to be considered. In all models, the moving solid    is described at time $t$ by $S(t) = h(t) + S$ for a fixed $S$. We assume that  $S$ has its lower tip at $r=0$, and  that in the vicinity of its lower tip, its surface is  described  by: 
$${   z = \gamma_S(r), \quad r \le r_0, \quad \theta \in (0,2\pi)      }$$ 
{for some  $r_0 \le 1$} and some Lipschitz function $\gamma_S$ with $\gamma_S(0) = 0$, $\gamma_S  \ge 0$. 
Notice that the solid velocity is given by $\dot{h}(t) \, e_z$.
  Similarly,  the  wall  $P$ is described in cartesian coordinates $(x,y,z)$, by 
$${  z = \gamma_P(x,y), \quad (x,y) \in \R^2,  } $$
for some Lipschitz function $\gamma_P$, with $\gamma_P(0,0) = 0$, $\gamma_P \le 0$.  
Accordingly, we denote the fluid domain 
$$ F(t) \: := \: \{\rmx = (x,y,z), \quad  \rmx \not\in \overline{S(t)}, \quad z > \gamma_P(x,y)\}.  $$
If $u=u(t,\rmx) = { \left(u_x(t,\rmx),u_y(t,\rmx),u_z(t,\rmx)\right)}$ and $p=p(t,\rmx)$ stand for the fluid velocity and pressure, the steady Stokes equations read 
\begin{equation} \label{Stokes}
 - \Delta u + \na p  = 0, \quad \div u = 0,  \quad  t> 0, \: \rmx \in F(t).
\end{equation}
We neglect gravity, as it plays no role in the discussion.  Our goal is to study the  force on the sphere,  that is 
{
\begin{equation} \label{defdrag}
{\cal F}_d(t) \: := \:  \int_{\pa S(t)} \left( 2D(u)n  - p n \right)d\sigma  \, \cdot \,  e_z.
\end{equation}
The notations $n$ and $D(u)$ refer to the normal vector pointing outside the fluid domain and the symmetric part of the
gradient respectively. 
}

\medskip
In order to determine ${\cal F}_d(t)$, one needs to specify the boundary conditions at the solid surface and at the plane. In all our  models for roughness, such conditions have the following general form:
\begin{equation} \label{BC1}
\bigl(u - \dot{h}(t) \, e_z \bigr) \cdot n\vert_{\pa S(t)} \: = \:  0,   \quad  \bigl(u - \dot{h}(t) \, e_z\bigr)  \times n \vert_{\pa S(t)} \: = \:  -  { 2 \beta_S \, [D(u) n]}   \times n\vert_{\pa S(t)} 
\end{equation}
and
\begin{equation} \label{BC2}
u  \cdot n\vert_{P} \: = \:  0,   \quad u  \times n \vert_{P} \: = \:  - { 2 \beta_P \,  [ D(u) n]}   \times n\vert_{P} 
\end{equation}
where $\beta_S, \beta_P \in [0,+\infty)$. These are boundary conditions of Navier type, the constants $\beta_P$ and $\beta_S$ being the slip lengths. Of special importance is the case $\beta_S = \beta_P = 0$, which corresponds to the no-slip condition.

\medskip
We  model the roughness in three different ways:
\begin{enumerate}
\item {\em through a lack of differentiability}. Namely, we consider a solid $S$ which is axisymmetric around $\{r=0\}$, and satisfies  
$$ \gamma_S(r) \: := \: 1 - \sqrt{1-r^2} + \eps r^{1+\alpha},  \quad \alpha \in [0,1), \quad r \le r_0. $$
This means that the solid surface is locally a smooth sphere $z =  1 - \sqrt{1-r^2}$, perturbed by a  less regular "rough profile"  of amplitude $\eps$.  For $\alpha = 0$, this profile is a spike, which has  Lipschitz regularity. For $\alpha > 0$, the profile is differentiable, with a H\"older derivative. 
For simplicity, we do not consider any roughness on the wall, and take the classical no-slip boundary conditions: 
 $\gamma_P = 0$, $\beta_P = \beta_S = 0$.

 \item {\em through a slip condition}. We consider the case of a ball $S$, of radius $1$,  falling vertically above a plane wall, with positive slip coefficients: 
$${ \gamma_S(r) = 1- \sqrt{1-r^2}}, \quad \gamma_P  = 0,  \quad  \beta_S, \beta_P > 0.$$
  Let us stress that such  modelling of the roughness by the addition of  (small) slip  is  commonly used. It is well-accepted in the context of rough hydrophobic surfaces \cite{Bocquet}, and a topic of debate in the context of hydrophilic ones {\it cf} \cite{Vinogradova2,Vinogradova1}. 
 \item {\em through a small parameter}.  Namely, the roughness is modelled through a small amplitude, high frequency perturbation of a plane wall. That means  $P $ is described by the equation
 $$ z = \gamma_P(x,y) \: :=  \: \eps \gamma(x/\eps,y/\eps), \quad  \eps \ll 1$$
 for some periodic and smooth  non-positive function $\gamma(X,Y)$, with $\gamma(0,0) = 0$. In parallel, we assume no roughness on the solid surface ($\gamma_S(r) = 1- \sqrt{1-r^2}$), and the classical no-slip conditions: $\beta_S = \beta_P = 0$.  
 \end{enumerate}
 
 \medskip
Note that if we take the parameters  $\eps$, $\beta_S$ and $\beta_P$  to be zero in the previous models, we are back to the classical situation of a curved and smooth solid falling towards a plane wall. The whole point is to derive the next order terms that are involved in the expression of ${\cal F}_d$. Note also that, in view of our models, the assumption that the solid translates along $r = 0$ is natural.  For the first two models, the whole geometry is axisymmetric. For the third one, one can consider  rough walls $P^\eps$ that are symmetric with respect to $x$ and $y$. In all these configurations, if the initial velocity field of the solid is along $r=0$,  both  the geometry and the Stokes flow inherit strong symmetry properties, forcing the  velocity field of the solid to be along  $r=0$ for all time. 

\medskip 
The ambition of this paper is to provide a rigorous and general methodology  to derive the drag term ${\cal F}_d$, in the 
 { regime}  of small distance $h$ between  the solid and the wall.  This methodology, which relies  on the calculus of variations,  will be explained in section \ref{methodology}. Then,  in section \ref{applications}, it will be applied  to our first two models  of rough surfaces. In this way, we will extend results from former formal computations, notably  those in \cite{Hocking, Feuillebois}. In  the last section \ref{oscillation}, we will turn to the third model of a small amplitude and high frequency boundary. This model is of particular interest, as it is connected to the phenomenon of apparent slip, which is a topic of current interest in fluid mechanics, see \cite{Lauga}. We will notably discuss the introduction of an effective slip length as a modelling for hydrophilic rough surfaces.

\section{Methodology for drag derivation} \label{methodology}
We present in this section a general approach to the derivation of  the drag force ${\cal F}_d(t)$  on the solid sphere $S(t)$. We first remark  that the geometric configuration at time $t$ is entirely characterized by the distance $h(t)$ between the lower tip of the solid and the origin $x = 0$. Thus, we can rewrite $ S(t) = S_{h(t)}, F(t) = F_{h(t)}$, with the family $(S_h, F_h)_h$ satisfying 
$$   S_h = h + S, \quad F_h = \Big\{ { \rmx, \quad \rmx} \not\in \overline{S_h}, \quad z > \gamma_P(x,y) \Bigr\}.$$ 
Moreover, considering the linear Stokes equation \eqref{Stokes} and boundary conditions \eqref{BC1}-\eqref{BC2}, we can write $ u(t,x)  =  \dot{h}(t)  u_{h(t)}(x)$ and $p(t,x)  =  \dot{h}(t) p_{h(t)}(x)$ where $u_h, p_h$ satisfy the steady problem 
\begin{equation} \label{Stokesbis}
 -\Delta u_h + \na p_h  = 0, \quad \div u_h = 0,  \quad  \rmx \in F_h 
\end{equation}
together with the boundary conditions 
\begin{equation} \label{BC1bis}
\bigl(u_h -  \, e_z \bigr) \cdot n\vert_{\pa S_h} \: = \:  0,   \quad   \bigl(u_h -  \, e_z\bigr)  \times n \vert_{\pa S_h} \: = \:   \, {  - 2 \beta_S [D(u_h) n]  }  \times n\vert_{\pa S_h} 
\end{equation}
and
\begin{equation} \label{BC2bis}
u_h  \cdot n\vert_{P} \: = \:  0,   \quad u_h  \times n \vert_{P} \: = \: { -2 \beta_P \,  [D(u_h)n] }  \times n\vert_{\pa P} 
\end{equation}
Accordingly, we can write 
$${\cal F}_d(t) = \dot{h}(t) {\cal F}_{h(t)}, \quad {\cal F}_h \: := \:  { \int_{\pa S_h} \left( 2D(u_h)n - p_h n \right)d\sigma  \, \cdot \,  e_z.} $$
The problem is to determine the behaviour of ${\cal F}_h$ in the limit $h \rightarrow 0$.  Our method to address this problem has three main steps: 
\begin{enumerate}
\item In a  first step,  we express the drag  ${\cal F}_h$ as the minimum of some energy functional.  One can do it using the variational interpretation of \eqref{Stokesbis}-\eqref{BC1bis}-\eqref{BC2bis}. It  allows to identify for  all our models of roughness  an energy functional  ${\cal E}_h$ and  a set of  "admissible fields" ${\cal A}_h$ such that 
$$ {\cal F}_h \: = \: \min_{u \in {\cal A}_h}  {\cal E}_h(u).  $$
The explicit definitions of ${\cal E}_h$  and ${\cal A}_h$ will be given at the end of this section. 
\item In a second step, we rely  on the minimization problem introduced in Step 1 to find an accurate lower bound for ${\cal F}_h$. Namely, we choose some appropriate  energy functional  $\tilde {\cal E}_h \le {\cal E}_h$ and some appropriate  set of admissible fields $\tilde {\cal A}_h \supset {\cal A}_h$ for which we can compute explicitly the minimimum and corresponding minimizer  $\tilde u$. In this way, we get  
$$ \tilde {\cal E}_h(\tilde u) \: = \: \min_{u \in \tilde {\cal A}_h}  \tilde {\cal E}_h(u) \: \le \: \min_{u\in {\cal A}_h}  {\cal E}_h(u)  =  
{\cal F}_h $$
which yields a lower bound. Of course, the relaxed functional  $\tilde {\cal E}_h$ and admissible set  $\tilde {\cal A}_h$ must remain close enough to the original ones, in order for this lower bound to be accurate.  We will  make them explicit for our various roughness models later on.  
\item In a third step, we choose some appropriate field $\check u \in {\cal A}_h$ so that 
$$ {\cal F}_h \: = \: \min_{u \in {\cal A}_h}  {\cal E}_h(u) \le   {\cal E}_h(\check u) $$
provides an accurate upper bound for the drag (that is with the same type of behaviour as the lower one). In many cases, as will be seen  later on, the minimizer  $\tilde u \in \tilde {\cal A}_h$  of the second step generally belongs   to the original set of admissible fields ${\cal A}_h$,  or at least can be slightly modified  to belong to ${\cal A}_h$. Thus, one can take in general $\check u \approx \tilde u$. 
\end{enumerate}

\medskip
Our goal in the present paper is to apply  this methodology to have a better understanding of roughness effects. 
In this section, we carry out step 1, that is the formulation of the drag in terms of some minimization problem. This step is very general, and  independent of the roughness issues. In the next sections, when turning to  step 2 and step 3, each roughness model will of course require specific calculations. 

\medskip
To link the drag to an extremum problem, we must distinguish between the case of  no slip ($\beta_S=\beta_P=0$) and the case of non-zero slip ($\beta_S > 0, \beta_P > 0$). 
\begin{itemize}
\item In the case of no-slip, { the divergence free-condition implies
$$
 \int_{F_h} | \na u_h |^2  = 2 \int_{F_h} |D(u_h)|^2.
$$}
{ Hence,} multiplying the Stokes  equation \eqref{Stokesbis}  by $u_h$ and integrating  over the fluid domain $F_h$, we obtain by Stokes formula
{
$$ 2 \int_{F_h} | D(u_h) |^2 \: = \: \int_{\pa S_h \cup P}  \left(  2D(u_h)n  - p_h n \right) \cdot   u_h  \, d\sigma\: = \: \int_{\pa S_h }  \left(  2D(u_h)n  - p_h n \right)d\sigma \cdot  e_z =  {\cal F}_h. $$
}
Moreover, we know that equation \eqref{Stokesbis} (together with the boundary conditions \eqref{BC1bis}-\eqref{BC2bis}) is the Euler equation of a minimization problem. Namely, 
$$ \int_{F_h} | \na u_h |^2  = \min \left\{  \int_{F_h} |\na u|^2,  \quad  u \in H^1_{loc}(F_h), \quad \na \cdot u = 0, \quad u\vert_P = 0, \: u\vert_{S_h} = e_z \right\} $$
(We remind that the Sobolev space $H^1_{loc}$  is the space of fields $u$ that are  locally square integrable, with distributional derivative $\na u$ also locally square integrable).
\footnote{{  As no abstract theory is needed in the remainder of the article, such mathematical details can be skipped without harm.}}   \\
Indeed, if $ u$ has the properties mentioned above, then $ u - u_h$ is zero along the boundary $\pa S_h \cup P$. So,   multiplying \eqref{Stokes} by $ u - u_h$  and integrating by parts, we end up with 
$$ \int_{F_h}  | \na u_h |^2  \: \le \: \int_{F_h}   \na u_h \, { : } \,   \na u \: \le \: \left( \int_{F_h}  | \na u_h |^2\right)^{1/2} \, \left( \int_{F_h}  | \na u|^2\right)^{1/2},$$
using the Cauchy-Schwarz inequality. The characterization of $u_h$ follows, and eventually yields that  ${\cal F}_h \: = \: \min_{u \in {\cal A}_h}   {\cal E}_h(u), $ with 
\begin{equation} \label{nrjnoslip}
{\cal E}_h(u) \: := \:  \int_{F_h} | \na u|^2, \quad {\cal A}_h\: := \:   \biggl\{u  \in H^1_{loc}(F_h), \quad  \na \cdot u = 0, \quad u\vert_P =0, \: u\vert_{\pa S_h} = e_z \biggr\}
\end{equation}

\medskip
\item   In the case of positive slip lengths $\beta_S, \beta_P$,  the computation is slightly different:
multiplying the Stokes equation by $u_h - u$, $ \: u \in {\cal A}_h$, we obtain after  integrating by parts: 
\begin{equation} \label{energynavier} 
\begin{aligned}
&  { \int_{F_h} 2 \,  D(u_h) : D(u_h - u)} + \: \frac{1}{\beta_S} \int_{\pa S_h}  \left( (u_h - e_z) \times n \right) \cdot   \left((u_h - u) \times n \right)  \, d\sigma \: \\
&  + \: \frac{1}{\beta_P} \int_{P}  \left( u_h  \times n \right) \cdot  \left( (u_h - u)  \times n \right)  \, d\sigma \: = \: 0.
\end{aligned}
\end{equation}
In order to recover full gradients instead of symmetric gradients, we proceed as follows. On one hand, by standard identities  of differential geometry (see for instance  \cite[Lemma 1, p. 233]{BrDeGe}), we have 
$$  D(v)n \times n = \frac{1}{2} \pa_n v \times n + \frac{1}{2}  v \times n  \: \mbox{ at $\pa S_h$ } $$
for any smooth  $v$   satisfying $v \cdot n = 0$ at $\pa S_h$. The last term at the r.h.s is connected  to the curvature of $\pa S$, which is simply $1$ by our choice of $S$. Similarly, 
$$  D(v)n \times n = \frac{1}{2} \pa_n v \times n  \: \mbox{ at $P$ } $$ 
for any smooth $v$ satisfying $v \cdot n = 0$ at $P$. On the other hand, writing 
$$ \Delta v = \div (\na v), \quad \mbox{respectively } \: \Delta v = 2 \div (D(v)), $$
 and integrating by parts, we get  that  for any smooth $v$  in ${\cal A}_h$ and any smooth $w$ satisfying $w \cdot n = 0$ at $\pa S_h \cup P$:
$$ \int_{F_h} \Delta v \cdot w =  - \int_{F_h}  \na v : \na w   \: + \:    \int_{\pa S_h} \left( \pa_n (v - e_z)   \times n \right)  \cdot (w   \times n)   \: + \:   \int_{P} \left( \pa_n v    \times n \right)  \cdot (w   \times n),  $$
respectively 
$$  \int_{F_h} \Delta v \cdot w  =  - \int_{F_h}  2 D(v) : D(w)   \: + \:    \int_{\pa S_h} \hspace{-0.1cm} \left( 2 D(v - e_z)n    \times n \right)  \cdot (w   \times n)   \: + \:   \int_{P} \left( 2 D(v)n     \times n \right)  \cdot (w   \times n). $$
{ Combining the previous identities, we get}
$$ \int_{F_h}  2 D(v) : D(w)  =  \int_{F_h}  \na v : \na w  \: +  \:    \int_{\pa S_h} \left( (v - e_z)   \times n \right)  \cdot (w   \times n). $$
We take $v = u_h$, $\: w = u - u_h$, and inject this last equality into \eqref{energynavier} to obtain
\begin{equation*}  
\begin{aligned}
&  \int_{F_h} |\na u_h|^2  \: + \: \left(\frac{1}{\beta_S} + 1\right) \int_{\pa S_h} | (u_h - e_z) \times n |^2    \, d\sigma \: + \: \frac{1}{\beta_P} \int_{P} | u_h \times n |^2   \, d\sigma  =  \int_{F_h} \na u_h : \na u  \\
+ & \: \left(\frac{1}{\beta_S} + 1\right) \int_{\pa S_h} \left(  (u_h - e_z) \times n \right) \cdot  \left(  (u - e_z) \times n \right) \, d\sigma \: + \: \frac{1}{\beta_P} \int_{P}  (u_h \times n) \cdot (u \times n) d\sigma. 
\end{aligned}
\end{equation*}
Use of the Cauchy-Schwartz inequality and of the Young inequality $\sqrt{ab} \: \le \: \frac{1}{2} (a + b)$ leaves us with 
\begin{equation*}  
\begin{aligned}
&  \int_{F_h} |\na u_h|^2  \: + \: \left(\frac{1}{\beta_S} + 1\right) \int_{\pa S_h} | (u_h - e_z) \times n |^2    \, d\sigma \: + \: \frac{1}{\beta_P} \int_{P} | u_h \times n |^2   \, d\sigma \\
\le &  \int_{F_h} |\na u|^2  \: + \: \left(\frac{1}{\beta_S} + 1\right) \int_{\pa S_h} | (u - e_z) \times n |^2    \, d\sigma \: + \: \frac{1}{\beta_P} \int_{P} | u \times n |^2   \, d\sigma .
\end{aligned}
\end{equation*}
Thus, we have  this time that ${\cal F}_h \: = \: \min_{u \in {\cal A}_h}   {\cal E}_h(u)$, with  
\begin{equation} \label{nrjslip}
\begin{aligned}
& {\cal E}_h(u) \: := \:  \int_{F_h} | \na u|^2 + { \: \left( \frac{1}{\beta_S} + 1\right)} \int_{\pa S_h} |(u - e_z) \times n |^2 \: + \: \frac{1}{\beta_P} \int_{P} |u \times n |^2, \\
& {\cal A}_h \: := \:  \biggl\{u \in H^1_{loc}(F_h), \quad   \na \cdot u = 0, \quad u \cdot n \vert_P =0, \: (u - e_z) \cdot n \vert_{\pa S_h} = 0 \biggr\}.
\end{aligned} 
\end{equation}
\end{itemize}
{ 
We note that contrary to the no-slip case, only the impermeability condition is included in the definition of the space ${\cal A}_h.$  
It can be shown that the Euler equation for the latter minimizing problem includes the boundary conditions \eqref{BC1bis}-\eqref{BC2bis} on the tangential part of the velocity-field 
by standard integration by parts as in the no-slip case.  {\em For brevity, we  shall replace the coefficient $1/\beta_S +1$ by $1/\beta_S$ in what follows}. This  means that  we shall include curvature effects in the slip coefficient.}
The characterization of the drag through energy functionals \eqref{nrjnoslip} and \eqref{nrjslip} will be applied to our first two roughness models in the next section.

\section{Application to various roughness models} \label{applications}
In this section, we detail the steps 2 and 3 of our methodology, both  in the case of a non-smooth boundary (model 1) and in the case  of slip boundary conditions (model 2). 

\subsection{ The case of non-smooth solids}
As emphasized in the introduction, we consider here the case of an axisymmetric solid $S$, whose boundary is described near its lower tip by 
$$ \gamma_S(r) = 1 - \sqrt{1-r^2} + \eps r^{1+\alpha}, \quad \alpha \in [0,1], \quad r \le r_0.     $$
The wall is flat, and no slip conditions are imposed at all boundaries. The drag is given by  
$${\cal F}_h \: = \: \min_{u \in {\cal A}_h}   {\cal E}_h(u), $$
 with the energy ${\cal E}_h$ and the set of admissible fields ${\cal A}_h$ given in \eqref{nrjnoslip}. 
 
\medskip
As the fluid domain $F_h$ is invariant by rotations around  $e_z$, much can be said about the minimizer   $u = u_h$. Indeed, for any rotation $R_\theta$ around $e_z$,  
$\: R_\theta u_h R_{-\theta}$ still belongs to  ${\cal A}_h$,  and has the same energy as $u_h$. Uniqueness of this minimizer yields
\begin{equation} \label{eq_invariance}
R_\theta \, u_h(R_{-\theta}\rmx) = u_h(\rmx), \quad \forall \, \rmx \in F_h.
\end{equation}
This means that $u_h$ has the following structure:
$$ u_h \: =  \: u_{h,r}(r,z) e_r \: + \: u_{h,\theta}(r,z) e_\theta \: + \: u_{h,z}(r,z) e_z, $$ 
where $(r,\theta,z)$, resp. $(e_r,e_\theta,e_z)$ are the cylindrical coordinates, resp. the cylindrical vector basis. One then remarks that $v_h = u_{h,r}(r,z) e_r \: + \: u_{h,z}(r,z) e_z$ still belongs to ${\cal A}_h$, with  ${\cal E}_h(v_h) \le  {\cal E}_h(u_h)$. Again, by uniqueness of the minimizer, we get $u_h = v_h$ and $u_{h,\theta} = 0$. Thus, the divergence free condition resumes to 
$$ \frac{1}{r} \pa_r (r u_{h,r}) + \pa_z u_{h,z} = 0.  $$
Together with  the boundary condition $u_{h,z}(r,0) = 0$, it leads to
\begin{equation} \label{stream}
 u_h = -\pa_z \phi \, e_r \: + \: \frac{1}{r}\pa_r (r \phi) \, e_z, 
 \end{equation}
with streamfunction $\phi(r,z) \: :=  \: -\int_0^z u_{h,r}(r,z')\, dz'$. 
The boundary conditions on $\phi$ are
\begin{equation} \label{BCstream}
\pa_z \phi\vert_{\pa S} = \pa_z \phi\vert_{P}  = 0,   \quad \pa_r(r \phi)\vert_{\pa S} = r,  \quad \phi\vert_{P} = 0.
\end{equation}
Thus, we can  without restriction include these last conditions in  the set of admissible fields: instead of the original definition in  \eqref{nrjslip},  we take
$$ {\cal A}_h\: := \:   \biggl\{u  \in H^1_{loc}(F_h), \quad u =  -\pa_z \phi \, e_r \: + \: \frac{1}{r}\pa_r (r \phi) \, e_z \:\:  \mbox{for some $\phi$ satisfying \eqref{BCstream}}\biggr\}. $$
We quote that the boundary conditions on $\phi$ at $\pa S$ yield
$$ \pa_z \phi(r,h+\gamma_S(r)) =0, \quad  \pa_r(r \phi)(r,h+\gamma_S(r)) = r, \quad r < r_0.$$
They imply in turn  that $ \phi(r,h+\gamma_S(r)) = \frac{r}{2} + \frac{c}{r} $ for some constant $c$. As $\phi(r,z) = -\int_0^z u_r(r,z')\,dz'$ is regular enough near $r=0$, we deduce $c=0$. Eventually
 \begin{equation} \label{BCstream2}
\begin{aligned} 
& \pa_z \phi(r,h+\gamma_S(r)) = 0,   \quad & \phi(r,h+\gamma_S(r)) = \frac{r}{2}, \quad r < r_0,\\
& \pa_z \phi(r,0) = 0,   \quad & \phi(r,0) = 0, \quad r < r_0.
\end{aligned}
\end{equation}

\medskip
From there, we obtain an accurate lower bound as follows.  Noticing that 
$$
|\nabla u|^2 = |\partial_{rz} \phi^2| + |\partial_z \phi/r|^2 + |\partial_{zz} \phi|^2 + |\partial_r[\partial_r(r\phi)/r]|^2 + |\partial_{rz}(r\phi)/r|^2.
$$ 
we anticipate  that in the limit of small $h$,  most of  the energy  ${\cal E}_h$  will come from a neighborhood of the lower tip of the sphere
$F^0_h := \{r < r_0 , \: 0 < z < h + \gamma_S(r)\},$
and will be due to the $z$ derivatives of the stream function $\phi$. 
  Accordingly,  we introduce the following relaxed minimizing set and energy functional:
 \begin{equation}
 \begin{aligned}
 \tilde{\cal A}_h \: & := \:   \biggl\{u  \in H^1(F^0_h), \quad u =  -\pa_z \phi e_r \: + \: \frac{1}{r}\pa_r (r \phi) e_z \:\:  \mbox{for some $\phi$ satisfying \eqref{BCstream2}}\biggr\}, \\
  \tilde{\cal E}_h \: & := \: \int_{F^0_h}  |\pa_z u_r|^2 =  \int_{F^0_h}  |\pa^2_z \phi|^2. 
  \end{aligned}
  \end{equation}

From the Euler equation $\pa^4_z \phi = 0$ and the boundary conditions  \eqref{BCstream2},  it follows easily that the latter minimum is realized with
$$
\tilde{\phi}_h(r,z) = \frac{r}{2} \Phi\left( \frac{z}{h + \gamma_S(r)}\right), \quad \text{ where } \quad \Phi(t) = t^2 (3-2t), \quad \forall \, t \in [0,1]
$$
and has for value:
\begin{equation} \label{eq_min}
\tilde{\cal F}_h = 6 \pi \int_0^{r_0} \frac{r^3 dr}{(h+ \gamma_S(r))^3}.
\end{equation}
We emphasize that this formula is general for no-slip boundary conditions. It does not require any special assumption on the solid surface.
In the case $\gamma_S (r) = 1 - \sqrt{1 - r^2} + \eps r^{1+\alpha},$ our lower bound satisfies
\begin{align*}
\tilde{\cal F}_h & \: = \: 6\pi \int_0^{r_0} \frac{r^3 dr}{(h+ \frac{r^2}{2} + \eps \,  r^{1+\alpha} + O(r^4))^3}  \\
 & \: = \: \frac{6\pi}{h} \,   {\cal I}\left(\eps h^{\frac{\alpha-1}{2}} \right) \: + \: O({\cal J}(\eps h^{\frac{\alpha-1}{2}},h)) + O(1)
 \end{align*}
 where 
\begin{equation} \label{eq_J}
{\cal I}\left(\beta\right) \: := \: \int_0^{\infty} \frac{s^3 ds}{(1+\frac{s^2}{2} +  \beta s^{1+\alpha})^3}, 
 \quad
 {\cal J}(\beta,h) \: := \: \int_0^{r_0/\sqrt{h}} \frac{s^7 ds}{(1+\frac{s^2}{2} +  \beta s^{1+\alpha})^4}.  
 \end{equation}
The computation of the asymptotic behaviours of $\cal I$ and $\cal J$  is detailed in  {\bf Appendix \ref{app_integrales}}. It yields the following results:
\begin{itemize}
 \item When $\beta  \ll 1$, we obtain :
\begin{equation}  \label{eq_kappabeta<<1}
 {\cal I}\left(\beta\right)\: = \: \frac{1}{1 + \lambda_{\alpha} \beta} + O(\beta^2), \quad 
 {\cal J}(\beta,h) \:= \:O(|\ln(h)|)
 \end{equation}
 with an explicit constant $\lambda_{\alpha}$ given in the appendix. 
\item When  $\beta \gg 1,$ { we have:}
\begin{equation} \label{eq_kappabeta>>1}
 {\cal I}\left(\beta\right) = 
\left\{
\begin{array}{lll}
{\mu_{\alpha}  \beta^{-\frac{4}{1+\alpha}}} &+ \: O\left(\dfrac{1}{\beta^3}\right), &   \text{for $\alpha >1/3$}, \\[10pt]
\dfrac{9}{4}\dfrac{\ln(\beta)}{\beta^{3}}  &+ \: O\left(\dfrac{1}{\beta^3}\right), & \text{for $\alpha =1/3,$} \\[10pt]
{ \mu_{\alpha}  \beta^{-\frac{2}{1-\alpha}}}  &+ \: O\left(\dfrac{1}{\beta^3}\right), & \text{for $\alpha < 1/3,$}
\end{array}
\right.
\end{equation}
where the value of $\mu_{\alpha}$ is also provided in the appendix. 
As regards the remainder, we have the following bound:
$$
{\cal J}\left(\beta,h\right) = O(|\ln(\beta) + \frac{1-\alpha}{2} \ln(h)|)
$$
\end{itemize}
Back to the drag force,   \eqref{eq_kappabeta<<1} and \eqref{eq_kappabeta>>1} yield  the following lower bound:  for  $\beta = \eps h^{\frac{\alpha-1}{2}} \ll 1$
\begin{equation} \label{eq_kappabeta<<<1}
\tilde{\cal F}_h =  \frac{6\pi}{h +\lambda_\alpha \eps h^{\frac{\alpha+1}{2}}}  \left( 1+ O\left(\beta\right)\right) \: + \:  O(|\ln(h)|)     
\end{equation}
and for $\beta = \eps h^{\frac{\alpha-1}{2}} \gg 1$
\begin{equation} \label{eq_kappabeta>>>1}
 \tilde{\cal F}_h =  
\left\{
\begin{array}{ll}
\dfrac{6 \pi \mu_{\alpha}}{ \, \eps^{\frac{4}{1+\alpha}} \, h^{\frac{3\alpha-1}{\alpha+1}}} \left( 1 + \beta^{\frac{1-3\alpha}{1+\alpha}} \right)
 + O(|\ln(h)|), & \text{ for $\alpha >1/3,$} \\[12pt] &\\
\dfrac{9 \pi  |\ln(h)|}{2 \, \eps^3} + O\left(\dfrac{|\ln(\eps)|}{\eps^3}\right), & \text{for $\alpha =1/3,$} \\[12pt] & \\
\dfrac{6\pi \mu_{\alpha} }{ \eps^{\frac{2}{1-\alpha}}}  \left( 1 + \beta^{\frac{3\alpha-1}{1+\alpha}} \right) + O(|\ln(\eps)|), & \text{for $\alpha <1/3,$} 
\end{array}
\right.
\end{equation}
Note that the expression  given in the case $\alpha < 1/3$ only matters  when $\eps \ll 1$ (otherwise, one  can just retain that $\tilde{\cal F}_h = O(1)$). 

\medskip
This concludes our study of a lower bound for the drag. Such bound is accurate, as we can with minor modifications obtain a similar upper bound. Indeed, it is an easy exercise to find a regular stream function  $\check{\phi}_h= \check{\phi}_h(r,z)$ defined on  $F_h$, equalling $\tilde \phi_h$ say on $F^0_h \cap \{r < r_0/2\}$, such that
$$\check{u}_h = \nabla \times (\check{\phi}_h e_{\theta}) \in {\cal A}_h, $$ 
and such that
$$ \int_{F_h} | \pa^2_z \check{\phi}_h |^2  = \int_{F^0_h} | \pa^2_z \tilde{\phi}_h |^2  \: + \: O(1) $$
uniformly in $h$ and $\eps$. We quote that the remainder term is uniformly bounded, because no singularity is created outside of the contact zone (that is outside a vicinity of $r = 0$).   We refer the reader to \cite{Hillairet07,HDGV10,HTT10} for more details in the case $\eps =0.$ Hence, we have: 
$$
\begin{array}{rcl}
{\cal F}_h &  \le  & {\cal E}_h(\check{u}_h)   \\
			& = &  \tilde{\cal F}_h + \displaystyle{\int_{F_h}} \left[  |\partial_{rz} \check{\phi}_h|^2 + |\partial_z \check{\phi}_h/r|^2 + |\partial_r[\partial_r(r\check{\phi}_h)/r]|^2 + |\partial_{rz}(r\check{\phi}_h)/r|^2\right] \: + \: O(1)\\[4pt]
			& = & \tilde{\cal F}_h + \displaystyle{\int_{F^0_h}} \left[  |\partial_{rz} \tilde\phi_h|^2 + |\partial_z \tilde\phi_h/r|^2 + |\partial_r[\partial_r(r\tilde\phi_h)/r]|^2 + |\partial_{rz}(r\tilde\phi_h)/r|^2\right] \: + \: O(1)
\end{array}
$$
The computation of the integral terms at the r.h.s. follows the lines of {\bf Appendix 
\ref{app_integrales}}. It yields some $O(\min(|\ln(h)|,|\ln(\eps)|))$ error term. The main reason for these integrals to be lower order terms is that in  the "curved" contact zone, the typical lengthscales in $z$ and $r$ are respectively $h$ and $\sqrt{h}$, so that $z$-derivatives are more singular than $r$-derivatives.  Eventually, 
$$ {\cal F}_h \: = \:  \tilde{\cal F}_h \: + \: O\left(\min(|\ln(\eps)|, |\ln(h)|)\right) $$
We stress that  this modelling of the roughness solves the famous  no-collision paradox discussed in the introduction. Indeed, as $h$ goes to zero for a given $\eps$, $\beta$ goes to infinity and the roughness effect yields a drag force ${\cal F}_h$ which is  always bounded by $c_\eps \, h^{-\gamma}$ for some $\gamma < 1$. In particular, it is weaker than in the smooth case, and the solid dynamics, which is governed by the o.d.e.  
$$  \ddot{h} \: + \:  \dot{h} \, {\cal F}_h   \: = \:  0$$
  allows for $h$ to cancel in finite time.

\subsection{The case of slip boundary conditions} 
We turn in this paragraph to our second model, in which roughness is involved through slip coefficients. We  want to have a close approximation of ${\cal F}_h =  \min_{{\cal A}_h} {\cal E}_h$, where this time ${\cal A}_h$ and ${\cal E}_h$ are defined in \eqref{nrjslip}.  We still have a rotational invariance in this case, so that we can again reduce 
${\cal A}_h$ by restricting to velocity fields of the type 
$$
u \: = \:  -\partial_z\phi \, e_r \: + \: \frac{1}{r} \partial_r [r\phi]  \, e_z. 
$$
The impermeability condition at $P$ yields again 
\begin{equation} \label{eq_BCi1}
\phi(r,0) = 0, \quad \forall \, r  \in (0,r_0).
\end{equation}
As regards $\pa S$,  we have for all $r \in (0,r_0)$
\begin{equation} \label{enun}
\begin{aligned}
 e_z \cdot n & = \dfrac{1}{\sqrt{1+ |\gamma'_S(r)|^2}} , \\
     u \cdot n & =  \dfrac{1}{\sqrt{1+ |\gamma'_S(r)|^2}}  \left( \gamma'_S(r) \pa_z\phi(r,h+\gamma_S(r)) +  \dfrac{1}{r} \pa_r (r \phi)(s,h+\gamma_S(r)) \right), \\
         &  = \dfrac{1}{\sqrt{1+ |\gamma'_S(r)|^2}} \dfrac{1}{r} \dfrac{\textrm{d}}{\textrm{d}r} [ r \phi(r,h+\gamma_S(r))] ,
 \end{aligned}
 \end{equation}
so that the impermeability condition leads to 
\begin{equation} \label{eq_BCi2}
\phi(r,h+\gamma_S(r)) = \dfrac{r}{2}, \quad \forall \, r \in (0,r_0).
\end{equation}
Accordingly, we introduce the relaxed set 
$$
\tilde{\cal A}_h := \left\{ u  \in H^1(F^0_h), \quad u =  -\partial_z\phi \, e_r \: + \: \frac{1}{r} \partial_r [r\phi]  \, e_z, \:  \mbox{$\phi$ satisfying } \:   \eqref{eq_BCi1} \: \mbox{ and  } \: \eqref{eq_BCi2} \right\}.
$$
with $F^0_h$ defined in the previous section.
We then need to define the relaxed energy  $\tilde{\cal E}_h$. As in the previous section, we  shall keep only $\partial_{zz} \phi$ in the gradient terms. But we shall not change the boundary integrals involved in \eqref{nrjslip}. Therefore, we compute:
\begin{itemize}
\item on $P$,  $\: u \times n = \partial_z \phi\, e_{\theta}$ 
\item on $\pa S \cap \{r < r_0\}, \: $ because of \eqref{enun}-\eqref{eq_BCi2}, 
$\: (u - e_z) \times n = \sqrt{1+ |\gamma'_S(r)|^2} \partial_z \phi \, e_{\theta}.  $
\end{itemize}
Hence, we introduce the approximate energy
\begin{multline*}
\tilde{\cal E}_h \: :=  \: \int_{F^0_h \cap \{r < r_0\}} |\partial_{zz} \phi (r,z)|^2 rdrd\theta dz \\
 + 2\pi \int_0^{r_0}  \left[ \dfrac{(1+|\gamma'_S(r)|^2)^{\frac{3}{2}}}{\beta_S}|\partial_z\phi(r,h+\gamma_S(r))|^2 + \dfrac{1}{\beta_P} |\partial_z \phi(r,0)|^2\right] rdr.
\end{multline*}
The corresponding mimimization problem is easy, because it amounts to find, for each value of $r < r_0$,  the minimizer of the functional 
$$ \tilde{\cal E}_h(r) = \int_0^{h+\gamma_S(r)} |\phi''_r(z)|^2  dz \\
 +  \left[ \dfrac{(1+|\gamma'_S(r)|^2)^{\frac{3}{2}}}{\beta_S}|\phi'_r(h+\gamma_S(r))|^2 + \dfrac{1}{\beta_P} |\phi'_r(0)|^2\right] $$
over functions $\phi_r = \phi_r(z)$ satisfying the inhomogeneous Dirichlet conditions 
$$ \phi_r(0) = 0, \quad  \phi_r(h+\gamma_S(r)) = \dfrac{r}{2} .  $$
This is a one-dimensional minimization problem, with Euler equation $\phi^{(4)}_r = 0$, endowed with above Dirichlet conditions, plus Robin type condition on  $\phi'_r$:
\begin{align*}
&  \phi''_r(h+\gamma_S(r)) + \dfrac{(1+|\gamma'_S(r)|^2)^{\frac{3}{2}}}{\beta_S}\phi'_r(h+\gamma_S(r))  = 0, \\
 & \phi''_r(0) -  \dfrac{1}{\beta_P} \phi'_r(0) = 0. 
 \end{align*}
After a few computations, the minimum of $\tilde {\cal E}_h$ is obtained for 
$$
\tilde{\phi}_h(r,z) = \dfrac{r}{2} \Phi\left(r, \dfrac{z}{h+\gamma_S(r)}\right),
$$
where $\Phi(r,t)$ is the  polynomial of degree $3$ in $t$ given by 
\begin{equation} \label{eq_chis}
\begin{aligned}
\Phi(r,t) := & - \dfrac{ 2 \  (\tbt + \tbt \ \tbb + \tbb )}{12 + 4 \ ( \tbt +\tbb ) + \tbb \ \tbt} \ t^3 + \dfrac{ 3 \ (2+ \tbt) \ \tbb }{12 + 4 \ ( \tbt +\tbb ) + \tbb \ \tbt} \ t^2 \\ 
& + \dfrac{ 6 \ (2+ \tbt) }{12 + 4 \ ( \tbt +\tbb ) + \tbb \ \tbt} \ t,
\end{aligned}
\end{equation}
where:
$$
\tbt = \tbt(r) \: := \:  \dfrac{(1+ |\gamma'_S(r)|^2)^{\frac{3}{2}} \, (h+ \gamma_S(r))}{\beta_S} , \qquad \tbb = \tbb(r) \: := \:  \dfrac{(h+ \gamma_S(r))}{\beta_P}.
$$
Note that the coefficients of $\Phi$ are uniformly bounded in $\alpha_S, \alpha_P$, that is  in $r < r_0$, $\beta_P,\beta_S,h$. In the limiting case $\beta_S=\beta_P=0$ (no-slip limit), we obtain formally $\Phi(t) = -2t^3 + 3 t^2$, in agreement with the computations of the previous section. 

\medskip
We now turn to the lower bound  
\begin{eqnarray*}
\tilde{\cal F}_h & = &  \min_{\tilde{\cal A}_h} \tilde{\cal E}_h \\
			& = & 2\pi \displaystyle{\int_0^{r_0}} \left[ \displaystyle{\int_0^1} |\partial_{tt} \Phi(r,s)|^2 ds + \tbt |\partial_t \Phi(r,1)|^2 + \tbb |\partial_t \Phi(r,0)|^2  \right] \dfrac{r^3 dr}{(h+ \gamma_S(r))^3}  	 		  	 
\end{eqnarray*}
We make the last  integral more  explicit by  replacing $\Phi$ by its value. We obtain 
$$
\tilde{\cal F}_h  \: =  \:  \frac{\pi}{2} \displaystyle{\int_0^{r_0}}  \left(I_1(r) + I_2(r)\right) \,    \dfrac{r^3 dr}{(h+ \gamma_S(r))^3}                     
$$
where the integrands $I_1$ and $I_2$  are given by 
\begin{align*}
I_1  \:  := \: &   \frac{12 \ \left(\tbt^2 \ \tbb^2 \: + \: 5 \  (\tbt^2 \ \tbb   +  \tbb^2 \ \tbt) \: + \:  
 4 \ (\tbt^2 + \tbb^2) \: + \: 20 \, \tbt \, \tbb \right)}{\left(12 \ + \  4 \ (\tbt \ + \  \tbb) \ + \ \tbt \tbb\right)^2}  \\
I_2  \:  := \: &   \frac{144 \, (\tbt \ + \ \tbb)}{\left(12 \ + \  4 \ (\tbt \ + \  \tbb) \ + \ \tbt \tbb\right)^2}\end{align*}
Note that $I_1$ and $I_2$ are uniformly bounded in $\alpha_S, \alpha_P$, that is  in $r < r_0$, $\beta_P,\beta_S,h$. Thus, expanding $\gamma_S,$ we obtain:
\begin{equation} \label{FI1I2}
\begin{array}{llcl} 
\tilde{\cal F}_h  & \: =  \: \dfrac{\pi}{2} \displaystyle{\int_0^{r_0}}  \left(I_1(r) + I_2(r)\right) \,    \dfrac{r^3 dr}{(h+ \frac{r^2}{2})^3} &+& O({\cal J}(0,h))   \\[6pt]
  & \: =  \:  \dfrac{\pi}{2} \displaystyle{\int_0^{r_0}}  \left(I_1(r) + I_2(r)\right) \,    \dfrac{r^3 dr}{(h+ \frac{r^2}{2})^3} &+& O(|\ln(h)|) 
   \end{array}                 
\end{equation}
where ${\cal J}(0,h)$ was introduced in \eqref{eq_J} and shown to be $O(|\ln(h)|).$
We must now distinguish between two cases, depending on the behaviour  of $h/\beta_S$ and  $h/\beta_P$:  
\begin{enumerate}
\item Either $h/\beta_S$ or $h /\beta_P$ is of order  $1$ or larger. Then, either $\alpha_S$ or $\alpha_P$ is of order $1$ or larger. It follows that  
$$ c \: \le \:  I_1(r) \: + \: I_2(r)  \: \le \: C, $$
for all $r < r_0$, where the constants $c,C$ are uniform with respect to all parameters. We then deduce from \eqref{FI1I2} that 
$$ \frac{c'}{h} \: \le \:  \tilde{\cal F}_h \: \le \: \frac{C'}{h}. $$
Note that in the limiting case $\beta_S=\beta_P=0$ (no-slip limit), we obtain formally:
 $$\alpha_S = \alpha_P = +\infty, \quad I_1 = 12, \quad I_2 = 0,  \quad \tilde{\cal F}_h = \frac{6\pi}{h}  $$ 
 recovering the classical result.  We also emphasize that the regime  considered here includes the case where  one of the slip coefficients  is zero. In particular, the drag force is stronger than  $c'/h$ in such a case, preventing any collision.  
 \item Both $h/\beta_S$ and $h /\beta_P$ are small.  This case requires more care. We first  notice that  
 \begin{equation} \label{DLalpha}
  \alpha_P = \frac{1}{\beta_P} \left(h + \frac{r^2}{2} + O(r^4) \right), \quad  \alpha_S = \frac{1}{\beta_S} \left(h + \frac{r^2}{2} + O(r^2(h +r^2)) \right). 
  \end{equation}
From there, for $r_0$ and $h$ small enough,  we get  
 $$ c \, J_1(r)\: \le \: I_1(r)  \: \le \:   C \, J_1(r), \quad J_1(r) \: := \: \left( \frac{a_P}{1+a_P} \: + \:  \frac{a_S}{1+a_S} \right)^2$$ 
 where $c,C > 0$ and  
 \begin{equation} \label{aSaP}
 a_P(r) \: := \:   \frac{1}{\beta_P} \left(h + \frac{r^2}{2}\right), \quad a_S(r) \: := \:   \frac{1}{\beta_S} \left(h + \frac{r^2}{2} \right).
\end{equation}
Then,  with the change of variable $r = \sqrt{h} u$ , we write
\begin{align*}
&   \displaystyle{\int_0^{r_0}} J_1(r) \,    \dfrac{r^3 dr}{(h+ \frac{r^2}{2})^3} \:\:  =  \: \\
&  \frac{1}{h} \int_0^{\frac{1}{2\sqrt{h}}}   \left( \frac{h/\beta_P\left(1+u^2/2\right)}{1 +  h/\beta_P\left(1+u^2/2\right)} \: + \: \frac{h/\beta_S\left(1+u^2/2\right)}{1 +  h/\beta_S\left(1+u^2/2\right)} \right)^2
\dfrac{u^3 du}{(1+u^2/2)^3}. 
\end{align*}
 In the regime of small $h/\beta_P$ and $h/\beta_S$,  we get that  this last integral is $o(1/\beta_P + 1/\beta_S)$. Finally, all of this  leads to
\begin{equation} \label{contribI1}  \frac{\pi}{2} \displaystyle{\int_0^{r_0}} I_1(r) \,    \dfrac{r^3 dr}{(h+ \frac{r^2}{2})^3}   \: = \:  o(1/\beta_P + 1/\beta_S).
\end{equation}

\medskip
It now remains to evaluate the contribution of $I_2$, which will yield the leading behaviour of $\tilde {\cal F}_h$. The use of \eqref{DLalpha} gives first
\begin{multline*}
 \frac{\pi}{2} \displaystyle{\int_0^{r_0}} I_2(r) \,    \dfrac{r^3 dr}{(h+ \frac{r^2}{2})^3}  \: = \:  \frac{\pi}{2}  \displaystyle{\int_0^{r_0}}   \frac{144 \, (a_S \ + \ a_P)}{\left(12 \ + \  4 \ (a_S \ + \  a_P) \ + \ a_S a_P\right)^2} \,    \dfrac{r^3 dr}{(h+ \frac{r^2}{2})^3}  \\ 
 + \:  O(1/\beta_P + 1/\beta_S).
 \end{multline*}
Then, straightforward manipulations show  that 
\begin{multline*}   \left( \frac{a_S}{(1+ c_1\, a_S)^2} \: + \:   \frac{a_P}{(1+ c_1 \, a_P)^2}  \right)   \:  + \:  O(J_1(r))     \: \le \: \frac{144 \, (a_S \ + \ a_P)}{\left(12 \ + \  4 \ (a_S \ + \  a_P) \ + \ a_S a_P\right)^2} \\ \le \:   
\left( \frac{a_S}{(1+ c_2 \, a_S)^2} \: + \:   \frac{a_P}{(1+ c_2 \, a_P)^2}  \right)   \:  + \:  O(J_1(r)) 
\end{multline*}
for some $c_1,c_2 > 0$. As seen in the treatment of $I_1$,  the $O(J_1(r))$ term will only contribute to the drag through a $o(1/\beta_P + 1/\beta_S)$ term.  The main contribution of $I_2$ to the drag will be governed by 
\begin{align*}
 & \frac{\pi}{2} \displaystyle{\int_0^{r_0}} \left( \frac{a_S}{(1+ c \, a_S)^2} \: + \:   \frac{a_P}{(1+ c \, a_P)^2}  \right)  \,    \dfrac{r^3 dr}{(h+ \frac{r^2}{2})^3}  \\
 &  = \:   \frac{\pi}{2 h} \int_0^{\frac{1}{2\sqrt{h}}}   \left( \frac{h/\beta_S \left(1+u^2/2\right)}{(1+ c \, h/\beta_S \left(1+u^2/2\right))^2} \: + \:   \frac{h/\beta_P \left(1+u^2/2\right)}{(1+ c \, h/\beta_P \left(1+u^2/2\right))^2}  \right)  \,    \dfrac{u^3 du}{(1 + \frac{u^2}{2})^3} \\
 & = \:   {\pi} \left\{ \frac{1}{\beta_S} \int_1^{1+\frac{1}{4h}} \frac{(x-1)dx}{(1+ c \, h/\beta_S \, x)^2\, x^2} \: + \:  \frac{1}{\beta_P} \int_1^{1+\frac{1}{4h}} \frac{(x-1)dx}{(1+ c \, h/\beta_P \, x )^2\, x^2}\right\}\\
 &  =   \:   {\pi} \left\{ \frac{1}{\beta_S} \int_1^{1+\frac{1}{4h}} \frac{dx}{(1+ c \, h/\beta_S \, x)^2\, x} \: + \:  \frac{1}{\beta_P} \int_1^{1+\frac{1}{4h}} \frac{dx}{(1+ c \, h/\beta_P \, x )^2 \, x}\right\} + O(1/\beta_P + 1/\beta_S)\\
 &  = \:  \pi \left( \frac{1}{\beta_S}+\frac{1}{\beta_P} \right) |\ln(h)| \: + \:   O(1/\beta_P + 1/\beta_S)
 \end{align*}
through standard manipulations. It yields  eventually  
\begin{equation} \label{contribI2}
  \frac{\pi}{2} \displaystyle{\int_0^{r_0}} I_2(r) \,    \dfrac{r^3 dr}{(h+ \frac{r^2}{2})^3}  \: = \: 
 \pi \left( \frac{1}{\beta_S}+\frac{1}{\beta_P} \right) |\ln(h)| \: + \:   O(1/\beta_P + 1/\beta_S)
 \end{equation}

 \bigskip
Combining \eqref{FI1I2}, \eqref{contribI1} and \eqref{contribI2}, we end up with the following lower bound for the drag:
$$ \tilde {\cal F}_h \: = \:  \pi \left( \frac{1}{\beta_S}+\frac{1}{\beta_P} \right) |\ln(h)| +  O(1/\beta_P + 1/\beta_S) + O(|\ln(h)|).  $$
\end{enumerate}

This lower bound is similar to the one derived by L.M. Hocking (see \cite{Hocking}).

\medskip
This concludes our study of a lower bound for the drag.  Hence, it remains to  obtain a similar upper bound. One could develop the same approach as in the previous section. Namely, one could look for  some  suitable extension $\check \phi_h$ of $\tilde{\phi}_h$, with similar behaviour for its energy.   However, due to the elaborate expression \eqref{eq_chis}, this would lead to tedious computations. We overcome this technical difficulty as follows:
\begin{enumerate}
\item  When $h/\beta_P$ or $h/\beta_S$ is of order 1 or larger, we take $\check u_h = \nabla \times (\check{\phi}_h e_{\theta})$, with the "no-slip" streamfunction $\check{\phi}_h$ built in the previous section. We obtain with this choice some $O(1/h)$ upper bound as expected. 
\item When $h/\beta_P$ and $h/\beta_S$ are small, a good  way to recover the right asymptotic behaviour  is to set 
$\check{u}_h \: := \:   \nabla \times [\check{\phi}_h(r,z) e_{\theta}]$  with
$$
\check{\phi}_h(r,z) \: := \:  \dfrac{r}{2} \Phi\left( r, \dfrac{z}{h+ \gamma_S(r)}\right),
$$
in $F_h^0$, where
$$
\Phi(r,t) \: := \:  \left( \dfrac{1}{1 + \alpha_P} + \dfrac{1}{1 + \alpha_S}  \right) \dfrac{t}{2} 
			     + \left( \dfrac{ \alpha_P}{1 + \alpha_P} + \dfrac{ \alpha_S}{1 + \alpha_S} \right) \dfrac{t^2}{2}.
$$
We extend then $\check{\phi}_h$ to the whole of $F_h$ with a stream function having bounded gradients. Calculations similar to the previous ones yield: 
$$ {\cal E}_h(\check{u}_h) \: = \: \pi \left( \frac{1}{\beta_S}+\frac{1}{\beta_P} \right) |\ln(h)| +  O(1/\beta_P + 1/\beta_S) \: + \:  O(|\ln h|) $$
where we insist that the $O(|\ln h|)$ is uniform with respect to $\beta_P$ and $\beta_S$. In particular, in the realistic regime of small slip lenghts, we obtain the exact same leading  behaviour  for the lower and upper bounds. For the sake of brevity, we leave the details to the reader. 
\end{enumerate}

\section{The case of a corrugated wall} \label{oscillation}
In this section, we focus on the third model of roughness described in the introduction, in which the wall has a small amplitude and high frequency oscillation: namely, 
$$  z \, =  \, \eps \gamma\left(\frac{x}{\eps}, \frac{y}{\eps}\right), \quad \eps > 0, \quad \gamma = \gamma(X,Y) \: \mbox{ 1-periodic,  } \: \gamma \le 0, \: \max \gamma = \gamma(0,0) = 0.  $$ 
We remind that the solid is assumed to be smooth,  and that no-slip conditions hold both at the solid surface and the wall. 
We shall pay special attention to  the regime  $\eps \ll h \ll 1$, that is when  the distance between the solid and the wall is much greater than the size of the roughness.

\medskip 
Such roughness model with a small parameter  is very popular, as it allows for multiscale analysis. This analysis has been notably performed in the context of wall laws. In this context, the idea is to replace the rough boundary by a flat one, and to impose there some good homogenized boundary condition, that expresses the mean effect of roughness.  This homogenization problem has been considered by physicists since  the early 90's, through numerics and explicit calculations for special geometries: see for instance \cite{Luchini:1995}. It has been adressed later on in some mathematical works, based on homogenization theory. We refer to \cite{Achdou:1998,Jager:2003} for periodic patterns of roughness, and to  \cite{DGV1,DGV2} for  random roughness.  The conclusion of these works is that, for small enough $\eps$,  one can replace the oscillating boundary by the flat one $\{z =0\}$, and impose there some Navier-type boundary condition: 
$$ u_z  = 0, \quad (u_x, u_y) \: = \: \eps B \, \pa_z (u_x, u_y)  $$ 
for some two by two positive matrix $B $, which is  sometimes called the "mobility tensor".  There has been a recent interest on qualitative properties of this tensor, for instance for shape optimization in microfluidics, {\it cf} \cite{Bazant}. 

\medskip
Another frequent idea is that a slip condition amounts to a no-slip condition at a shifted wall. Combining this idea with the previous one, some recent articles have  suggested a drag force of the type ${\cal F}_h  \sim \frac{1}{h+ \beta \eps}$ for some positive $\beta$: see \cite{Vinogradova2,Feuillebois}. We will discuss  this result in a rigorous manner here.

\medskip
First,  one can use  the methodology of  section \ref{methodology} to derive some lower and upper bounds.  

\begin{itemize}
\item As regards the lower bound, let us show that 
 $${\cal F}_h  \: \ge \:   \frac{6\pi}{h+ \lambda \eps} \: + \:  O(|\ln(h + \lambda \eps)|), \quad \mbox{ for } \: \lambda \: := \:  -\min \gamma > 0.$$
Indeed, we  have 
 $$ {\cal F}_h = \min_{u \in {\cal A}_h}  {\cal E}_h(u)  $$ 
where ${\cal E}_h(u)$ and  ${\cal A}_h$ are given by \eqref{nrjnoslip}. Let us now define 
$$P^\lambda  \: := \:  \{ z = -\eps \lambda\}, \quad F_h^\lambda \: :=  \: \{ \rmx,  \quad \rmx  \not\in \overline{S}_h, \quad z > -\eps \lambda \}. $$
Any field $u$ of ${\cal A}_h$ can be extended  by zero below the rough wall  so that it can be seen as an element of the larger set
$$ \tilde {\cal A}_h \: :=  \:   \biggl\{u  \in H^1_{loc}(F_h^\lambda), \quad  \na \cdot u = 0, \quad u\vert_{P^\lambda} =0, \: u\vert_{\pa S_h} = e_z \biggr\}. $$
Then, obviously, 
$$ {\cal F}_h  \: \ge \: \min_{u \in \tilde {\cal A}_h}  {\cal E}_h(u).  $$
But the r.h.s of this inequality is exactly the drag force associated to  the (smooth)  solid $S_h$ and the (smooth) plane 
$P^\lambda$. As the distance between the two is $h+\lambda \eps$, we deduce  the expected lower bound. 
\item  With similar arguments, one has the upper bound: 
$$  {\cal F}_h  \: \le \:   \frac{6\pi}{h} \: + \:  O(|\ln(h)|).$$
Indeed, let us define this time 
$$P^0  \: := \:  \{ z = 0\}, \quad F_h^0 \: :=  \: \{ \rmx,  \quad \rmx  \not\in \overline{S}_h, \quad z > 0 \}.$$  
Let $u^0_h$ be  the Stokes flow in the smooth domain $F_h^0$. As the distance between $S_h$ and $P^0$  is $h$, the drag force satisfies 
$$ \int_{F_h^0} |\na u^0_h|^2 \: = \:  \frac{6\pi}{h} \: + \:  O(|\ln(h)|) $$
Now, $u^0_h$ can be extended by zero below $P^0$ and defines in this way an element of ${\cal A}_h$. In particular, 
$ {\cal F}_h  \: \le \:  \int_{F_h} |\na u^0_h|^2.$
\end{itemize}
Hence, our methodology allows to derive quickly the  inequalities 
$$  \frac{6\pi}{h+ \lambda \eps}  \: + \:  O(|\ln(h + \lambda \eps)|) \: \le \:  {\cal F}_h \: \le \:  \frac{6\pi}{h} \: + \:  O(|\ln(h)|).$$
Interestingly, these bounds are satisfied for any r\'egime of parameters $\eps$ and $h$. In particular, it provides the right asymptotic when $\eps$ and $h$ are of the same order.

\medskip
Nevertheless, when $\eps \ll h$, it is fair to notice that a multiscale analysis gives a much refined description of the drag. For the sake of completeness, we briefly present it here. It relies on an asymptotic expansion of the Stokes flow  $u_h=u^\eps_h$ with respect to $\eps$. This expansion has already been described in close contexts, for instance in \cite{Jager:2003}, and we only recall its main elements. To keep track of the $\eps$ dependency, we write $P^\eps$ instead of $P$, $F^\eps_h$  instead of $F_h$. We denote again $P^0$ and $F^0_h$ their smooth counterparts. 
 
\medskip
 The basic  idea is  to build an approximate solution $u^\eps_{h,app}(\rmx)$ of  \eqref{Stokesbis}-\eqref{BC1bis}-\eqref{BC2bis},  in the form of an expansion in powers of $\eps$: 
 \begin{equation} \label{expansion}
u^\eps_{h,app}(\rmx) \: = \:   u^0_h(\rmx) \: + \: \eps \biggl( u^1_h(\rmx) + U^1_h(x,y, \rmx/\eps) \biggr) \: + \:  \dots \: + \: \eps^N \biggl( u^N_h(\rmx) + U^N_h(x,y, \rmx/\eps) \biggr) 
\end{equation}
Each term of this expansion  has two parts: 
\begin{itemize}
\item A regular part $u^i_h = u^i_h(\rmx)$ which models the macroscopic variations of the solution.
\item A boundary layer correction  $U^i_h = U^i_h(x,y, \rmX) $, which accounts for the fast variations of the solution near the oscillating boundary. Hence, it depends on the macroscopic variables $x,y$, but also on the microscopic variable $\rmX = \rmx /\eps$.  It is defined for all
$$x,y \in \R^2, \quad \rmX=(X,Y,Z) \: \mbox{ such that } \:    Z >  \gamma(X,Y). $$
Moreover,  $U^i_h$ is periodic in $X,Y$ (due to the periodicity of the rough bondary $\gamma$ in $X,Y$)  and satisfies 
$$ \lim U^i_h(x,y, X,Y,Z) = 0, \quad \mbox{as } \:  Z \rightarrow +\infty $$
Back to the original variable $\rmx$, this last condition corresponds to a boundary layer of typical size $\eps$ near the rough wall $P$. 
\end{itemize}
Accordingly,  the corresponding  pressure field should read
$$ p^\eps_{h,app}(\rmx) \: =  \: p^0_h(\rmx) + P^0_h(x,y,\rmx/\eps)  +  \eps \biggl( p^1_h(\rmx) + P^1_h(x,y,\rmx/\eps) \biggr) \: + \: \dots $$ 

\medskip
We  remind here the derivation of the $O(1)$ and $O(\eps)$ terms, which are enough for our purpose. First, if we inject the above expansions in \eqref{Stokesbis}-\eqref{BC1bis}-\eqref{BC2bis}, and let $Z \rightarrow +\infty$,we obtain  by standard manipulations that 
$$ -\Delta u^0_h + \na p^0_h  = 0 \: \mbox{ and } \:  \div u^0_h = 0 \: \mbox{  in } \:  F^0_h, \quad u^0_h\vert_{\pa S_h} = e_z, \quad  u^0_h\vert_{P^0} = 0.  $$
Thus, we recover as expected that  the leading term  of the expansion is the Stokes flow without roughness. We then extend $u^0_h$  and $p^0_h$ by zero below $P^0$, so that they are defined over the whole $F^\eps_h$. Such extensions trivially satisfy the Stokes equation for $ z < 0$, as well as the no-slip condition at $P^\eps$. Moreover, the velocity is continuous across the plane $\{z = 0\}$. But there is a $O(1)$ jump in the normal derivative. This explains the introduction of a boundary layer corrector with amplitude $O(\eps)$. Indeed, its  gradient has  amplitude $O(1)$, and allows to correct  this artificial jump. 

\medskip
Let us introduce the following notations: 
\begin{align*}
V_h(x,y,\rmX) \: & :=  \: u^1_h(x,y,0) + U^1_h(x,y,\rmX), \quad  Z > \gamma(X,Y), \\
P_h(x,y,\rmX) \: & := \: P^0_h(x,y,\rmX) - p^0_h(x,y,0), \quad  Z > 0,  \\
P_h(x,y,\rmX) \: & := \: P^0_h(x,y,\rmX), \quad  0 > Z > \gamma(X,Y).
\end{align*}

Note  that, following  the expansion \eqref{expansion}, 
$ u^0(x,y,0) + \eps V_h(x,y,\rmx/\eps) $ should be an approximation of the whole flow in the boundary layer. 
Plugging \eqref{expansion} into the equations, we derive formally the following Stokes system:
\begin{equation*}
\left\{
\begin{aligned}
-\Delta_{\rmX} V_h + \na_{\rmX} P_h \:  =  \: 0,  &  \quad Z > \gamma(X,Y), \: Z \neq 0, \\ 
 \na_{\rmX} \cdot  V_h \:  =  \: 0, & \quad    Z > 0 , \: Z \neq 0, \\  
 V_h \:  = \: 0, &  \quad   Z = \gamma(X,Y). 
 \end{aligned}
 \right. 
 \end{equation*}
 together with the jump conditions 
 $$ V_h\vert_{Z=0^+} - V_h\vert_{Z=0^-} = 0, \qquad \left(\pa_Z V_h - P_h e_Z\right)\vert_{Z=0^+} -  \left(\pa_Z V_h - P_h e_Z\right)\vert_{Z=0^-} = - \pa_z u^0_h(x,y,0). $$
 Again, we stress that these jump conditions ensure the smoothness of the  whole flow across the artificial boundary $\{ z = 0\}$. Note  that by the divergence-free condition $\pa_z u^0_{h,z}(x,y,0) = 0$, so that only the horizontal components of $\pa_z u^0_h(x,y,0)$  are non-zero. Let us  also point out that the variables $x,y$ are just parameters in the system. In other words, one has 
$$V_h(x,y,\rmX) = {\cal V}(\rmX)  \pa_z u^0_h(x,y,0), \quad P_h(x,y,\rmX) = {\cal P}(\rmX) \cdot  \pa_z u^0_h(x,y,0)$$
  for some 3-by-3  matrix function ${\cal V}$ and some 3d vector ${\cal P}$ which satisfy  the (matricial) Stokes system
\begin{equation} \label{BL}
\left\{\begin{aligned}
-\Delta_{\rmX} {\cal V} + \na_{\rmX} {\cal P} \:  =  \: 0,  &  \quad Z > \gamma(X,Y), \\ 
 \na_{\rmX} \cdot  {\cal V} \:  =  \: 0, & \quad     Z > \gamma(X,Y), \\ 
 {\cal V} \:  = \:  0,  \, &  \quad   Z = \gamma(X,Y), 
 \end{aligned}
 \right.
 \end{equation}
 together with 
 \begin{equation*}
 {\cal V}\vert_{Z=0^+} - {\cal V}\vert_{Z=0^-} = 0, \qquad \left(\pa_Z {\cal V}- {\cal P} \otimes e_Z
 \right)\vert_{Z=0^+} -\left(\pa_Z {\cal V}- {\cal P} \otimes e_Z
 \right)\vert_{Z=0^-}    = - \left( \begin{smallmatrix} 1 & 0 & 0 \\ 0 & 1 & 0 \\ 0 & 0 & 0 \end{smallmatrix} \right) 
 \end{equation*}
This system of pde's, depending only  on $\rmX$, with periodic boundary conditions in the horizontal variable $X,Y$, has been extensively studied.  We remind the following proposition, extracted from \cite{Jager:2003}:
\begin{proposition}
The solution ${\cal V}$ of system \eqref{BL}  converges exponentially at infinity, that is 
$$ \left| {\cal V}(X,Y,Z) - {\cal V}^\infty \right| \: \le \: C \, e^{-\delta Z} $$
for some constant 3-by-3 matrix ${\cal V}^\infty$ and some $\delta > 0$. Moreover,  ${\cal V}^\infty$ is of the form 
$$ {\cal V}^\infty \: = \: \begin{pmatrix} {\cal B} & \begin{matrix} 0 \\ 0 \end{matrix} \\ \begin{matrix}0 & 0 \end{matrix} & 0 \end{pmatrix}  $$ 
for some symmetric positive definite 2-by-2 matrix ${\cal B}$. 
\end{proposition}
The non-zero block ${\cal B}$ is sometimes called the mobility tensor, see \cite{Bazant}. We   stress that ${\cal B}$ is symmetric and definite positive, so diagonalizable in an orthonormal basis with positive eigenvalues. This fact will be used below.

\medskip
Back to  $V_h$, we obtain that 
$$ \bigl( V_{h,x}, V_{h,y} \bigr)(x,y,\rmX) \: \rightarrow \: {\cal B} \, \pa_z\left(u^0_{h,x}, u^0_{h,y}\right)(x,y,0), \quad V_{h,z}(x,y,\rmX) \rightarrow 0, \quad \mbox{as } \: Z \rightarrow +\infty. $$  
As the boundary layer correction $U^1_h$ should decay at infinity, we obtain the boundary condition for the macroscopic correction $u^1_h$ at $P^0$. That is 
$$ (u^1_{h,x}, u^1_{h,y})  = {\cal B} \pa_z\left(u^0_{h,x}, u^0_{h,y}\right), \quad u^1_{h,z} = 0, \quad \mbox{ at } P^0. $$
Together with the Stokes equations
$$-\Delta u^1_h + \na p^1_h  = 0 \: \mbox{ and } \:  \div u^1_h = 0 \: \mbox{  in } \:  F^0_h, $$
and the boundary condition at the solid surface
$$ u^1_h\vert_{\pa S_h} = 0 $$  
this determines $u^1_h$, and ends the derivation of the $O(\eps)$ term of the expansion. The next order terms solve the same kind of equations, with inhomogeneous data coming from lower order profile.   

\medskip
In a second step, one can show rigorously that the approximate solution $u^\eps_{h,app}$ is close to  the exact solution $u^\eps_h$. Indeed, introducing the differences   
$$v \: :=  \: u^\eps_{h,app} - u^\eps_h, \quad  \mbox{ and  } \quad  q \: :=  \: p^\eps_{h,app} - p^\eps_h$$
leads to 
$$  - \Delta v + \na q = R^\eps_h, \quad \div v = r^\eps_h, \quad v\vert_{\pa S_h}  = \varphi^\eps_h. $$
with remainder terms $R^\eps_h$, $r^\eps_h$ and $\varphi^\eps_h$. For instance, the boundary data $\varphi^\eps_h$ is due to the  boundary layer terms $U^i(\rmx,\rmx/\eps)$,  that  do not vanish at $\pa S_h$.   We stress that the assumption  $\eps \ll h$ is crucial for these remainders to be small.  First,  the boundary layer corrections decay exponentially over a typical lengthscale $\eps$. To  make it exponentially small at $\pa S_h$, one needs $\eps \ll h $.  Moreover, all other remainder terms are small with respect to $\eps$,  but  diverging with respect to $h$. Very roughly, they behave like $O((\eps/h)^N)$ where $N$ is the number of terms in the expansion \eqref{expansion}. The diverging powers of $h$ come from taking derivatives of the $u^i_h$, which are singular with respect to $h$. 
Again, the smallness condition $\eps \ll h$ is necessary.

\medskip
 From there, as the remainder terms are small, one can through energy estimates deduce the smallness of $v$, that is $u^\eps_h \approx u^\eps_{h,app}$. In particular, for $\eps$ small enough compared to $h$, the drag force on $\pa S_h$ reads 
$$ {\cal F}_h \: = \:  \int_{\pa S_h} \left( \frac{\pa u^\eps_h}{\pa n } - p^\eps_h n \right) \cdot e_z =   \int_{\pa S_h} \left( \frac{\pa (u^0_h + \eps u^1_h}{\pa n } - (p^0_h + \eps p^1_h)  n \right) \cdot e_z \: + \: o(\eps). $$
Moreover, it is easily seen that  the fields ${\check u}^\eps_h \: := \: u^0_h + \eps u^1_h$, $\: {\check p}^\eps_h \: := \: p^0_h + \eps p^1_h$ satisfy the Stokes equation in $F_h^0$, together with the boundary conditions
$$ {\check u}^\eps_h\vert_{\pa S_h} = e_z, \quad \mbox{ and } \quad {\check u}^\eps_{h,z} \vert_{P^0}  = 0, \quad  
\left( {\check u}^\eps_{h,x},  {\check u}^\eps_{h,y}\right)\vert_{P^0} \: = \eps {\cal B} \, \pa_z\left( u^0_{h,x},  u^0_{h,y}\right)\vert_{P^0}. $$
By the axisymmetry of $F^0_h$,    $R_\theta \, u^0_h \: = \:  u^0_h R_\theta$ for any horizontal rotation $R_\theta$. As ${\cal B}$  is symmetric definite positive, this allows us to assume, up to a change of orthonormal basis, that ${\cal B}$ is diagonal with positive coefficients $\beta_x, \beta_y$. Now, there are two ways to interpret the effect of roughness. 
\begin{itemize}
\item On one hand, one can write the latter boundary condition as a slip condition of Navier type:
$$ \left( {\check u}^\eps_{h,x}, {\check u}^\eps_{h,y}\right)\vert_{P^0} \: = \eps {\cal B} \, \pa_z\left( {\check u}^\eps_{h,x},  {\check u}^\eps_{h,y}\right)\vert_{P^0} \: +  \: o(\eps).$$
This is the so-called phenomenon of apparent slip, see \cite{Lauga}. In the isotropic case $\beta := \beta_x = \beta_y$, one can use the bounds  on the drag force derived in section \ref{applications}. In the r\'egime $\eps^{-1} h \gg 1,$ this yields: 
\begin{equation} \label{bound1}
\frac{c}{h}   \le {\cal F}_h \le \frac{C}{h}, \qquad c, C > 0.
\end{equation}
\item On the other hand,  the drag force reads 
\begin{align*}
 {\cal F}_h \: & \approx \:  \int_{\pa S_h} \left( \frac{\pa u^0_h}{\pa n } - p^0_h n \right) \cdot e_z + \eps  \int_{\pa S_h} \left( \frac{\pa u^1_h}{\pa n } - p^1_h n \right) \cdot e_z \\
& = \:   \int_{\pa S_h} \left( \frac{\pa u^0_h}{\pa n } - p^0_h n \right) \cdot e_z  + \eps  \int_{F^0_h} \na u^1_h : \na u^0_h \\
& = \:  \int_{\pa S_h} \left( \frac{\pa u^0_h}{\pa n } - p^0_h n \right) \cdot e_z +  \eps  \int_{P^0} u^1_h \cdot
\left( \frac{\pa u^0_h}{\pa n } - p^0_h n \right)\\
& = \:  \int_{\pa S_h} \left( \frac{\pa u^0_h}{\pa n } - p^0_h n \right) \cdot e_z -  \eps  \int_{P^0} \left( \beta_x |\pa_z u^0_{h,x}|^2 + \beta_y |\pa_z u^0_{h,y}|^2 \right) 
\end{align*}
Using again the symmetry properties of  $u^0_h$, we obtain that 
\begin{equation} \label{dragapprox}
 {\cal F}_h \:  \approx \:  \int_{\pa S_h} \left( \frac{\pa u^0_h}{\pa n } - p^0_h n \right) \cdot e_z \: - \:  \eps \, \beta  \int_{P^0} \left(  |\pa_z u^0_{h,x}|^2 +  |\pa_z u^0_{h,y}|^2 \right), \quad \beta \: := \: \frac{\beta_x + \beta_y}{2}.   
 \end{equation}
But  this last expression can be seen  as the drag force created by a Stokes flow $u^{\beta}_h$, between the solid $S_h$ and a shifted wall $P_\beta \: := \: z  = -  \eps \, \beta$. Indeed, following \cite{Feuillebois}, we get
\begin{align*}
  \int_{\pa S_h} \left( \frac{\pa u^\beta_h}{\pa n } - p^\beta_h n \right) \cdot e_z & =   \int_{\pa S_h} \left( \frac{\pa u^\beta_h}{\pa n } - p^\beta_h n \right) \cdot u^0_h \\
& = \int_{F^0_h} \na u^\beta_h : \na u^0_h \\
& = \int_{S^0_h} e_z \cdot \left( \frac{\pa u^0_h}{\pa n } - p^0_h n \right)  + \int_{P^0} u^\beta_h \cdot    \left( \frac{\pa u^0_h}{\pa n } - p^0_h n \right) 
\end{align*}
Using that 
$$ u^\beta_h(x,y,0) = u^\beta_h(x,y,-\eps\beta) + \eps \beta \pa_z u^\beta(x,y,-\eps \beta) + o(\eps) = 
\eps \beta \pa_z u^0_h(x,y,0) + o(\eps)  $$ 
we recover the same expression as in \eqref{dragapprox}. This interpretation of the roughness  effect as a shift of the smooth  wall yields 
\begin{equation} \label{bound2}
 {\cal F}_h \approx \frac{6\pi}{h+\eps\beta}. 
 \end{equation}
\end{itemize}
Note that bounds \eqref{bound1} and  \eqref{bound2} are coherent, since we are here within the asymptotics  $\eps \ll h \ll 1$. 
Of course, these bounds would lead to very different behaviours for smaller $h$, as \eqref{bound1} forbids collision whereas \eqref{bound2} allows it.

\appendix

\section{Asymptotics of ${\cal I}$ and ${\cal J}$}
\label{app_integrales}
In this appendix, we detail the computation of ${\cal I}(\beta)$ and ${\cal J}(\beta,h)$ defined in \eqref{eq_J} depending on the values of $\beta.$

\paragraph{Case $\beta <<1.$} When $\beta$ is small, we expand with respect to $\beta$ :
$$
 \dfrac{s^3}{(1+\frac{s^2}{2}+\beta s^{1+\alpha})^3} = \dfrac{s^3}{(1+\frac{s^2}{2})^3} - 3 \beta \dfrac{s^{4+\alpha}}{(1+\frac{s^2}{2})^4} + O\left( \beta^2 \dfrac{s^{5+2\alpha}}{(1+\frac{s^2}{2})^5} \right)
$$
This yields:
$$
{\cal I}(\beta) = \int_0^{\infty} \dfrac{s^3 \, ds}{(1+\frac{s^2}{2})^3}  - \beta \int_0^{\infty} \dfrac{3 s^{4+\alpha} \, ds }{(1+\frac{s^2}{2})^4} + O\left(\beta^2\right)
$$
where routine calculations yield:
$$
 \int_0^{\infty} \dfrac{s^3 \, ds}{(1+\frac{s^2}{2})^3} = {1} , \quad \int_0^{\infty} \dfrac{3 s^{4+\alpha} \, ds }{(1+\frac{s^2}{2})^4} = \dfrac{2^{\frac{\alpha+1}{2}}\pi (3+\alpha)(1-\alpha^2)}{8 \cos\left( \frac{\pi \alpha}{2}\right)} =:\lambda_{\alpha}. 
$$
Replacing in ${\cal I}(\beta)$, we obtain, 
$$
{\cal I}(\beta) = {1} - {\lambda_{\alpha} \beta} + O(\beta^2) = \dfrac{1}{1 + \lambda_{\alpha} \beta} + O(\beta^2) 
$$

\paragraph{Case $\beta >>1.$} When $\beta$ is large, we split ${\cal I}(\beta) = {\cal I}^0(\beta) + {\cal I}^{\infty}(\beta)$ where:
$$
{\cal I}^0(\beta) = \int_0^{1}  \dfrac{s^3\,ds}{(1+\frac{s^2}{2}+\beta s^{1+\alpha})^3}, \quad 
{\cal I}^\infty_{\beta} = \int_1^{\infty}  \dfrac{s^3\,ds}{(1+\frac{s^2}{2}+\beta s^{1+\alpha})^3}, 
$$
To compute ${\cal I}^{\infty}(\beta),$ we set $s= \beta^{\frac{1}{1-\alpha}} \tilde{s}$ and expand the integrand with respect to $1/(\beta^{\frac{2}{1-\alpha}}).$
This yields:
$$
{\cal I}^{\infty}(\beta) = \dfrac{1}{\beta^{\frac{2}{1-\alpha}}} \int^{\infty}_{\beta^{-\frac{1}{1-\alpha}}} \dfrac{\tilde{s}^3 \, d\tilde{s}}{(\frac{\tilde{s}^2}{2} + \tilde{s}^{1+\alpha})^3} + O\left(\frac{1}{\beta^{4}} \right)$$
Consequently, we distinguish three cases :
\begin{itemize}
\item for $\alpha > 1/3,$ \: ${\cal I}^{\infty}(\beta) = O\left( \dfrac{1}{\beta^3}\right)$\\
\item for $\alpha =1/3,$ \: ${\cal I}^{\infty}(\beta) = \dfrac{3}{2} \dfrac{\ln (\beta)}{\beta^3}  \: + \: O\left(\dfrac{1}{\beta^3}\right) $ \\
\item for $\alpha < 1/3,$ \: 
$
{\cal I}^{\infty}(\beta) =   \dfrac{\mu_\alpha}{\beta^{\frac{2}{1-\alpha}}} \,    \: + \: O\left(\dfrac{1}{\beta^3}\right), \quad \mu_\alpha := \displaystyle{\int_0^{\infty} }\dfrac{ \tilde{s}^3 \, d\tilde{s}}{(\frac{\tilde{s}^2}{2} + \tilde{s}^{1+\alpha})^3}.
$
\end{itemize}
In ${\cal I}^0(\beta),$ we set $u = \beta s^{1+\alpha}$ and expand the integrand w.r.t. $1/\beta^{\frac{1}{1+\alpha}}.$ This yields:
{ $$
{\cal I}^0(\beta) = \dfrac{1}{(1+\alpha)\beta^{\frac{4}{1+\alpha}}} \, \int_0^{\beta} \dfrac{u^{\frac{3-\alpha}{1+\alpha}} \, du}{(1+u)^3} + O\left(\frac{1}{\beta^4}\right).
$$}
with :
\begin{itemize}
\item for $\alpha < 1/3,$ \:  ${\cal I}^0(\beta) = O\left(\dfrac{1}{\beta^3}\right)$\\
\item for $\alpha =1/3,$ \:  ${\cal I}^0(\beta) = \dfrac{3}{4} \dfrac{\ln (\beta)}{\beta^3} \: + \: O\left(\dfrac{1}{\beta^3}\right)$\\
\item for $\alpha > 1/3,$ \:  ${\cal I}^0(\beta) = \dfrac{\mu_\alpha}{\beta^{\frac{4}{1+\alpha}}}  \: + \: O\left(\dfrac{1}{\beta^3}\right), \quad \mu_\alpha:= \dfrac{1}{1+\alpha} \displaystyle{\int_0^{\infty}}  \dfrac{u^{\frac{3-\alpha}{1+\alpha}} \, du}{(1+u)^3} $
\end{itemize}
We obtain \eqref{eq_kappabeta>>1} comparing the values of ${\cal I}^0(\beta)$ and ${\cal I}^{\infty}(\beta)$ in the three cases $\alpha < 1/3,$ $\alpha=1/3$ and $\alpha > 1/3.$

It remains to handle the remainder term ${\cal J}(\beta,h).$ As previously,
we split it into ${\cal J}^0(\beta,h) + {\cal J}^{\infty}(\beta,h)$ where:
$$
{\cal J}^0(\beta,h) = \int_0^1 \dfrac{s^7ds}{(1+\frac{s^2}{2} + \beta s^{1+\alpha})^4} = O(1),
\quad
{\cal J}^{\infty}(\beta,h) = \int_{1}^{{r_0}/{\sqrt{h}}}  \dfrac{s^7ds}{(1+\frac{s^2}{2} + \beta s^{1+\alpha})^4}
$$ 
We set $s = \beta^{\frac{1}{1-\alpha}} \tilde{s}$ in the last integral. This yields:
$$
{\cal J}^{\infty}(\beta,h) \le  \int_{1/\beta^{1/(1-\alpha)}}^{r_0/(\sqrt{h} \beta^{1/(1-\alpha)})}  \dfrac{\tilde{s}^7ds}{(\frac{\tilde{s}^2}{2} +  \tilde{s}^{1+\alpha})^4}
= 
\left\{
\begin{array}{lcl}
O(|\ln(h)|) & \text{if $\beta \ll 1$}\\[6pt]
O(|\ln(\beta) + \frac{1-\alpha}{2}\ln h |) & \text{if $\beta \gg 1$}.
\end{array}
\right.
$$
\bibliographystyle{abbrv}
\bibliography{biblio-mod}
   
\end{document}